\documentclass[11pt,a4paper,reqno]{amsart}

\usepackage{amsmath}
\usepackage{amsfonts}
\usepackage{amssymb}
\usepackage{latexsym}
\usepackage{bbm}

\renewcommand\geq{\geqslant}
\renewcommand\leq{\leqslant}

\newcommand{\bea}{\begin{eqnarray*}}
\newcommand{\bean}{\begin{eqnarray}}
\newcommand{\eea}{\end{eqnarray*}}
\newcommand{\eean}{\end{eqnarray}}
\newtheorem{theorem}{Theorem}[section]
\newtheorem{lemma}[theorem]{Lemma}
\newtheorem{conjecture}[theorem]{Conjecture}

\newcommand{\A}{\mathcal{A}}
\newcommand{\B}{\mathcal{B}}
\newcommand{\F}{{\mathcal{F}}}
\renewcommand{\H}{{\mathcal{H}}}
\newcommand{\I}{\mathcal{I}}
\newcommand{\N}{\mathbbm{N}}

\newcommand{\sn}{S_n}

\newcommand{\Matroids}[3]{\mathbbm{M}_{#2}^{#1}{(#3)}}
\newcommand{\NIMatroids}[3]{\mathcal{M}_{#2}^{#1}{(#3)}}
\newcommand{\size}[1]{|{#1}|}

\newcommand{\binary}[1]{\vec{#1}}

\newcommand{\Stirling}[2]{\left\{ {#1} \atop {#2} \right\}}

\newcommand{\T}{\mathcal{T}}
\newcommand{\Refine}{\mbox{\scshape{Refine}}}
\newcommand{\Random}{\mbox{\scshape{Random}}}
\newcommand{\Expand}{\mbox{\scshape{Expand}}}

\newcommand{\Free}{\mbox{Free}}

\title{On the number of matroids on a finite set}
\author{W. M. B. Dukes}
\thanks{Supported by EC's Research Training Network `Algebraic Combinatorics in Europe',
grant HPRN-CT-2001-00272 while the author was at Universit\`{a} di Roma Tor Vergata, Italy.}
\email{dukes@labri.fr}
\address{LaBRI, Universit\'{e} Bordeaux 1, 351 cours de la Lib\'{e}ration, 33405 Talence Cedex, France.}
\begin{document}
\maketitle
\begin{abstract}
In this paper we highlight some enumerative results concerning matroids of low rank 
and prove the tail-ends of various sequences involving the number of matroids on a finite set to be log-convex.
We give a recursion for a new, slightly improved, lower bound on the number of 
rank-$r$ matroids on $n$ elements when $n=2^m-1$.
We also prove an adjacent result showing the point-lines-planes conjecture to be 
true if and only if it is true for a special sub-collection of matroids.
Two new tables are also presented, giving the number of paving matroids on at most eight elements.
\end{abstract}

\section{Definitions}
Let $\sn$ be a finite set of size $n$.
A matroid $M$ is a pair $M(\sn,\mathcal{B})$ such that $\mathcal{B}$ is a non-empty collection 
of subsets of $\sn$ satisfying the following basis-exchange property:
$$\forall\; X,Y \in\mathcal{B},\;\; \forall x\in X-Y, \exists y\in Y-X \mbox{ s.t } X-\{x\} \cup\{y\} \in\mathcal{B}.$$
$\sn$ is called the {\it{ground set}} of the matroid and $\B$ the collection of {\it{bases}}.
A consequence of this axiom is that all sets in $\mathcal{B}$ have the same cardinality, 
called the {\it{rank}} of the matroid.
The collection of {\it{independent}} sets of $M$ is ${\mathcal{I}}(M) = \{I\,:\,I\subseteq B, B\in\mathcal{B}(M)\}$.
An element $x\in\sn$ is called a {\it{loop}} if $\{x\}\not\in\mathcal{I}(M)$. 
If for every pair of distinct elements $\{x,y\}\subseteq \sn$
there exists $I\in\I(M)$ s.t. $I\supseteq \{x,y\}$, then $M$ is called a {\it{simple}}.

Let $\Matroids{k}{r}{\sn}$ be the class of rank-$r$ matroids on $\sn$ 
with all $k$-element sets independent, non-empty for all $0\leq k \leq r \leq n$,
and the union of these classes over $k\leq r\leq n$ is denoted $\Matroids{k}{}{\sn}$.
We say two matroids $M_1(\sn,\B_1)$ and $M_2(\sn,\B_2)$ are {\it{isomorphic}} if there 
exists a permutation $\pi$ of $\sn$ such that $\B_2 = \{ \pi(B) : B\in\B_1\}$.
Let $\NIMatroids{k}{r}{\sn}$ and $\NIMatroids{k}{}{\sn}$ be the corresponding classes of non-isomorphic matroids (i.e. `different' matroids).  
The collections $\Matroids{0}{r}{\sn}$, $\Matroids{1}{r}{\sn}$ and $\Matroids{2}{r}{\sn}$ 
identify the classes of rank-$r$ matroids, loopless matroids and simple matroids, respectively, on $\sn$.
$\Matroids{r-1}{r}{\sn}$ is the class of rank-$r$ paving matroids on $\sn$ and the single matroid in 
$\Matroids{r}{r}{\sn}$ is the uniform rank-$r$ matroid on $n$ elements.
Note that $\Matroids{i+1}{r}{\sn} \subset \Matroids{i}{r}{\sn}$ for all $0\leq i <r$.
Given $M\in\Matroids{0}{r}{\sn}$ and $X\subseteq \sn$, the {\it{restriction}} of $M$ to $X$, 
denoted $M|X$, is the matroid with independent sets $\{I \cap X : I\in\mathcal{I}(M)\}$.
The {\it{dual}} of $M$ is the matroid $M^{\star}\in\Matroids{0}{n-r}{\sn}$ 
with $\B(M^{\star}) := \{ \sn-B : B \in \B(M)\}$.

The rank of $A\subseteq \sn$ is defined by $r(A):= \max\{|X|: X\subseteq A, X\in\I(M)\}$.
We call $A\subseteq \sn$ a {\it{flat}} if $r(A\cup\{x\}) = r(A)+1$ for all $x\in\sn -A$.
We denote by $\F_k(M)$ the collection of rank-$k$ flats of $M$ and the numbers $W_k(M):= \size{\F_k(M)}$ are called the {\it{Whitney numbers of the second kind}}.
We may also define a matroid in terms of its flats: the pair $M(\sn,\F = \{\F_0,\F_1,\ldots ,\F_r\})$ is
a matroid if for all $0\leq i <r$ and $F_i \in\F_i$, the collection
$\{ X-F_i : X\in\F_{i+1}\mbox{ and }X\supset F_i\}$ partitions the set $\sn -\F_i$.

For a general matroid $M$, the relationship between the collection of flats and independent sets 
may be seen as follows: a set $F$ is in $\F(M)$ if $I\cup\{x\} \in\I(M)$ whenever 
$I\subseteq F$, $I\in\I(M)$ and $x\not\in F$. 
Conversely, a set $I$ is in $\I(M)$ if for every $x \in I$, $I-F =\{x\}$ for some $F\in\F(M)$.

The upper truncation $U_k(M)$ of $M\in\Matroids{0}{r}{\sn}$ is the matroid obtained 
by removing the rank-$i$ flats of $M$ for $k\leq i<r$. The resulting matroid is a rank-$k$ matroid.
Given $N \in\Matroids{0}{r}{\sn}$ and $M\in\Matroids{0}{r+1}{\sn}$, we say $M$ is an {\it{erection}} of $N$ if $U_{r}(M)=N$.

A collection $\H$ of subsets of $\sn$ is a $d$-partition of $\sn$ if every $d$-set of $\sn$ 
is contained in a unique set of $\H$ and $|H|\geq d$ for all $H\in\H$.
We call the $d$-partition {\it{trivial}} if $\H = \{\sn\}$.

\section{Enumerative results and conjectures}
Several papers and books addressing the number of matroids on a finite set have used the terms `number of matroids' 
where `number of non-isomorphic matroids' or `number of non-isomorphic simple matroids' was intended. 
The tightest bounds on the total number of matroids 
(i.e., fixing $n$ and summing over $r$) are to be found in Knuth~\cite{knuthpaving} and Piff and Welsh~\cite{piffwelsh}. 
Some surprising recent work gives the asymptotics for the number of binary matroids, see Wild~\cite{wild}.
Also, the long-standing conjecture of Welsh~\cite{welshconj} that 
$\size{\NIMatroids{0}{}{S_{m+n}}} \geq \size{\NIMatroids{0}{}{S_m}}\size{\NIMatroids{0}{}{S_n}}$
has recently been proved independently by Lemos~\cite{lemos} and Crapo and Schmitt~\cite{craposchmitt}.

To resolve the confusion regarding what the term `number of matroids' means, 
in this paper we will distinguish between the six classes: 
${\Matroids{0}{r}{\sn}}$,
${\Matroids{1}{r}{\sn}}$,
${\Matroids{2}{r}{\sn}}$,
${\NIMatroids{0}{r}{\sn}}$,
${\NIMatroids{1}{r}{\sn}}$ and
${\NIMatroids{2}{r}{\sn}}$. 

To align our notation with the the standard literature on this subject, let us mention the following;
Crapo~\cite{singleelement} gives the first 4 numbers in the 
$\size{\NIMatroids{2}{}{\sn}}$ row of the table for $\size{\NIMatroids{2}{r}{\sn}}$ in Section 5.
The table for $\size{\NIMatroids{2}{r}{\sn}}$ was first given in Blackburn, Crapo and Higgs~\cite{blackcrapohiggs} 
and is reprinted in Welsh~\cite[p.306]{welshbook}.
The table $\size{\NIMatroids{0}{r}{\sn}}$ was first given in Acketa~\cite{acketa84}, 
extending earlier results from Acketa~\cite{acketa79}, and is reprinted at the back of Oxley's book~\cite{oxleybook}. 
The numbers enumerated in Acketa~\cite{acketa78} were $\size{\NIMatroids{0}{2}{\sn}}$. 

Let $p(i)$ be the number of partitions of the integer $i$, $\Stirling{n}{i}$ 
the Stirling numbers of the second kind and $b(i)$ the $i^{\mbox{\tiny th}}$ Bell number.
Through consideration of the loops and points of general matroids, we have the following relations:
$\size{\NIMatroids{0}{0}{\sn}} = 1$, 
$\size{\Matroids{0}{0}{\sn}} = 1$, 
$\size{\Matroids{1}{1}{\sn}} = 1$,
$\size{\NIMatroids{1}{1}{\sn}} = 1$,
$\size{\Matroids{2}{n}{\sn}}=1$, 
$\size{\NIMatroids{2}{n}{\sn}}=1$, 

\begin{eqnarray*}
\size{\Matroids{0}{r}{\sn}} &=& \sum_{i=r}^n {n\choose i} \size{\Matroids{1}{r}{S_i}},\\
\size{\Matroids{1}{r}{\sn}} &=& \sum_{i=r}^{n} \Stirling{n}{i} \size{\Matroids{2}{r}{S_i}},\\
\size{\NIMatroids{0}{r}{\sn}} &=& \sum_{i=r}^{n} \size{\NIMatroids{1}{r}{S_i}}.
\end{eqnarray*}

Using these equalities, we have the following closed forms for certain low-rank matroids:

\begin{theorem}[Dukes~\cite{dukesphd}]\label{numbers}
For all $n\geq 2$,
\renewcommand{\theenumi}{\roman{enumi}}
\begin{enumerate}
\item $\size{\NIMatroids{0}{1}{\sn}} = n$, $\size{\NIMatroids{0}{2}{\sn}} = p(1)+\cdots +p(n)-n$,
\item $\size{\Matroids{0}{1}{\sn}} = 2^n -1$, $\size{\Matroids{0}{2}{\sn}} = b(n+1)-2^n$,
\item$\size{\Matroids{1}{2}{\sn}} = b(n)-1$, $\size{\NIMatroids{1}{2}{\sn}} = p(n)-1$.
\end{enumerate}
\renewcommand{\theenumi}{\arabic{enumi}}
\end{theorem}

By duality we can also enumerate high-rank matroids as the following lemma shows:

\begin{lemma}\label{tail-simple}
For all $n\geq 5$, 
\renewcommand{\theenumi}{\roman{enumi}}
\begin{enumerate}
\item $\size{\Matroids{2}{n-1}{\sn}}=2^{n}-1-{n+1\choose 2}$, 
\item $\size{\NIMatroids{2}{n-1}{\sn}}=n-2$,
\item $\size{\Matroids{2}{n-2}{\sn}}=b(n+1)+{n+3\choose 4}+2{n+1\choose 4}-2^n -{n+1\choose 2}2^{n-1}$,
\item $\size{\NIMatroids{2}{n-2}{\sn}}=p(1)+\ldots +p(n)+6-4n$.
\end{enumerate}
\renewcommand{\theenumi}{\arabic{enumi}}
\end{lemma}

\begin{proof}
A matroid  $M\in\Matroids{2}{n-1}{\sn}$ is specified by a unique circuit of length $k$ for some $3\leq k \leq n$, 
the elements not contained in this circuit are the isthmuses of $M$. 
If $\B(M^{\star}) =\{\{x_1\},\ldots ,\{x_k\}\}$ then this circuit is $C=\{x_1,\ldots ,x_k\}$ and
the number of ways to choose such a circuit is simply the number of subsets of $S_n$ of size at least 3, hence 
$\size{\Matroids{2}{n-1}{\sn}} = 2^n-1-n-{n\choose 2}$.
The corresponding non-isomorphic number is $\size{\NIMatroids{2}{n-1}{\sn}} = (n+1)-3=n-2$.

If $M\in\Matroids{0}{2}{\sn}$ then $M^{\star}\in\Matroids{0}{n-2}{\sn}$.
Suppose $\F_0(M)=\{F_0\}$ and $\F_1(M) = \{H_1,\ldots ,H_k\}$. 
Let $F_i = H_i - F_0$ for all $1\leq i \leq k$. 
Then $\{F_0,F_1,\ldots ,F_k\}$ partitions the set $\sn$. 
$M^{\star}$ is not simple iff there exists $X\in{\sn\choose n-2}$ such that for 
all $(f_i,f_j)\in F_i\times F_j$ with $0<i<j$, $\{f_i,f_j\}\not\subset X$.
It transpires that $M^{\star}$ is not simple if 
\begin{itemize}
\item $k=2$ and ($\size{F_1}=1$, $\size{F_2}\geq 1$) or ($\size{F_1}=2$, $\size{F_2} \geq 2$),
\item $k=3$, $\size{F_1}=\size{F_2}=1$ and $\size{F_3}\geq 1$,
\end{itemize}
and relabellings thereof.
Enumerating this class of matroids whose duals are not simple gives:
\begin{eqnarray*}
{\sum_{i=0}^{n-2} {n\choose i} (n-i) - {n\choose n-2} + 3{n\choose 4} + \sum_{i=5}^n {n\choose i}{i\choose 2}
+ {n\choose 3} + \sum_{i=4}^n {n\choose i}{i\choose 2}}\\
\;=\; {n+1 \choose 2} 2^{n-1} - {n+3 \choose 4} -2 {n+1\choose 4}.
\end{eqnarray*}
Thus 
\begin{eqnarray*}
\size{\Matroids{2}{n-2}{\sn}} 
&=& \size{\Matroids{0}{2}{\sn}} - \left\{{n+1 \choose 2} 2^{n-1} - {n+3 \choose 4} -2 {n+1\choose 4}\right\} \\
&=& b(n+1)-2^n -{n+1 \choose 2} 2^{n-1} +{n+3 \choose 4} +2{n+1\choose 4}
\end{eqnarray*}
and the corresponding non-isomorphic number is
\begin{eqnarray*}
\size{\NIMatroids{2}{n-2}{\sn}} &=& \size{\NIMatroids{0}{2}{\sn}} - (n-1)-(n-3)-(n-2) \\
&=& p(1)+\ldots +p(n)-n-3(n-2).
\end{eqnarray*}
\end{proof}

A well cited conjecture (proposed by Welsh~\cite{welshcombbook}) states that the number of `different' 
matroids on $\sn$ (i.e. $\size{\NIMatroids{0}{r}{\sn}}$) is unimodal in rank. 
It seems unfair to restrict this conjecture to only one of the six classes of matroids so we propose;
\begin{conjecture}
The sequences 
$\left\{\size{\Matroids{0}{r}{\sn}}\right\}_{r=0}^{n}$,
$\left\{\size{\Matroids{1}{r}{\sn}}\right\}_{r=1}^{n}$,
$\left\{\size{\Matroids{2}{r}{\sn}}\right\}_{r=2}^{n}$,
$\left\{\size{\NIMatroids{0}{r}{\sn}}\right\}_{r=0}^{n}$,
$\left\{\size{\NIMatroids{1}{r}{\sn}}\right\}_{r=1}^{n}$ and
$\left\{\size{\NIMatroids{2}{r}{\sn}}\right\}_{r=2}^{n}$
are unimodal and attain their maximum at $r=\lceil n/2 \rceil$.
\end{conjecture}
The following theorem summarizes partial progress with this conjecture.
\begin{theorem}[Dukes~\cite{dukes1}] For all $n\geq 6$,
\renewcommand{\theenumi}{\roman{enumi}}
\begin{enumerate}
\item $\size{\Matroids{0}{0}{\sn}}\leq \size{\Matroids{0}{1}{\sn}}\leq \size{\Matroids{0}{2}{\sn}}\leq 
	\size{\Matroids{0}{3}{\sn}}$,
\item $\size{\NIMatroids{0}{0}{\sn}}\leq \size{\NIMatroids{0}{1}{\sn}}\leq 
	\size{\NIMatroids{0}{2}{\sn}}\leq \size{\NIMatroids{0}{3}{\sn}}$,
\item $\size{\Matroids{1}{1}{\sn}}\leq \size{\Matroids{1}{2}{\sn}}\leq \size{\Matroids{1}{3}{\sn}}$,
\item $\size{\NIMatroids{1}{1}{\sn}}\leq \size{\NIMatroids{1}{2}{\sn}}\leq \size{\NIMatroids{1}{3}{\sn}}$,
\item $\size{\Matroids{2}{2}{\sn}}\leq \size{\Matroids{2}{3}{\sn}}$,
\item $\size{\NIMatroids{2}{2}{\sn}}\leq \size{\NIMatroids{2}{3}{\sn}}$.
\end{enumerate}
\renewcommand{\theenumi}{\arabic{enumi}}
\end{theorem}

In problems involving unimodal sequences, it is sometimes easier to show a sequence to 
be log-concave and therefore unimodal.  It is not clear that the sequences 
$\left\{\size{\Matroids{0}{r}{\sn}}\right\}_{r=0}^{n}$,
$\left\{\size{\Matroids{1}{r}{\sn}}\right\}_{r=1}^{n}$,
$\left\{\size{\Matroids{2}{r}{\sn}}\right\}_{r=2}^{n}$,
$\left\{\size{\NIMatroids{0}{r}{\sn}}\right\}_{r=0}^{n}$,
$\left\{\size{\NIMatroids{1}{r}{\sn}}\right\}_{r=1}^{n}$ and
$\left\{\size{\NIMatroids{2}{r}{\sn}}\right\}_{r=2}^{n}$
are log-concave. The evidence in the tables would seem to suggest log-concavity. 
However, this is only the case for $n$ small.
We exhibit the following theorem which shows the beginning (in some cases the end) of these sequences to be 
log-convex once $n$ is sufficiently large. First we require two lower bounds to aid the proof.

\begin{theorem}[Dukes~\cite{dukes2}]\label{slc-bound}
For all $2<r<n$, $\size{\Matroids{r-1}{r}{\sn}} > 2^{{n\choose r}/2n}$ 
and $\size{\NIMatroids{r-1}{r}{\sn}} > {1\over n!} 2^{{n\choose r}/2n}$.
\end{theorem}

\begin{theorem}
The matroid numbers satisfy the following inequalities
\renewcommand{\theenumi}{\roman{enumi}}
\begin{enumerate}
\item $\size{\Matroids{0}{0}{\sn}} \size{\Matroids{0}{2}{\sn}} \;>\; \size{\Matroids{0}{1}{\sn}}^2$ for all $n\geq 4$,
\item $\size{\NIMatroids{0}{0}{\sn}} \size{\NIMatroids{0}{2}{\sn}} \;>\; \size{\NIMatroids{0}{1}{\sn}}^2$ for all $n\geq 9$,
\item $\size{\Matroids{1}{1}{\sn}} \size{\Matroids{1}{3}{\sn}} \; >\; \size{\Matroids{1}{2}{\sn}}^2$ for all $n\geq 94$,
\item $\size{\NIMatroids{1}{1}{\sn}} \size{\NIMatroids{1}{3}{\sn}} \;>\; \size{\NIMatroids{1}{2}{\sn}}^2$ for all $n\geq 67$,
\item $\size{\Matroids{2}{n-2}{\sn}} \size{\Matroids{2}{n}{\sn}} \; >\; \size{\Matroids{2}{n-1}{\sn}}^2$ for all $n\geq 11$,
\item $\size{\NIMatroids{2}{n-2}{\sn}} \size{\NIMatroids{2}{n}{\sn}} \;>\; \size{\NIMatroids{2}{n-1}{\sn}}^2$ for all $n\geq 8$.
\end{enumerate}
\renewcommand{\theenumi}{\arabic{enumi}}
\end{theorem}

\begin{proof}
Each follows from combining known closed forms, and in some cases inequalities, for the numbers 
and using induction.
The first and second inequalities follow directly from Theorem~\ref{numbers}.
The third inequality follows since
$\size{\Matroids{1}{1}{\sn}}\size{\Matroids{1}{3}{\sn}} >
\size{\Matroids{2}{3}{\sn}} > 2^{(n-1)(n-2)/12}$, by using Theorem~\ref{slc-bound} with $r=3$.
It can then be verified that $2^{(n-1)(n-2)/12} > (b(n)-1)^2 = \size{\Matroids{1}{2}{\sn}}^2$ for all $n\geq 94$.
Similarly, the fourth inequality follows from applying the same method 
but by using the non-isomorphic lower bound of Theorem~\ref{slc-bound}.
The final two inequalities follow from Lemma~\ref{tail-simple}.
Some of the inequalities also hold for smaller values of $n$.
\end{proof}

The above sequences are much more likely to resemble the sequence $\{2^{n\choose r}\}_{r=0}^{n}$ 
which is log-convex at both ends of the sequence and has a log-concave `cap' in the middle.

\section{Lower bound on the number of paving matroids}
The tightest asymptotic lower bound for the number of matroids was given in Knuth~\cite{knuthpaving}.
The method involved proving a lower bound on the number of rank-$\lfloor n/2\rfloor$ paving matroids on $\sn$.
Theorem~\ref{slc-bound} showed his lower bound inequality to hold for $2<r<n$ (~\cite{knuthpaving} 
proved the case $r=\lfloor n/2 \rfloor$)
giving a lower bound on the number of $(r-1)$-partitions of $\sn$. 
The class of non-trivial $(r-1)$-partitions corresponds directly to the class of rank-$r$ paving matroids.
It was noted in Dukes~\cite{dukes2} that experimental results suggested a particular class of 
the partitions to be larger than each of the other classes. 
The following lemma gives a recursion for the enumeration of this class, 
hence the tightest lower bound possible using this approach when the size of the ground set is of the form $2^m-1$.

\begin{theorem}
If $n=2^m-1$ for some $m\in\N$, 
then $\size{\Matroids{r-1}{r}{\sn}} \geq 2^{\size{\mathcal{U}(r,n)}}$ for all $r \geq 4$
where $\size{\mathcal{U}(r,n)}$ satisfies the recursion:
\begin{eqnarray*}
\size{\mathcal{U}(r+1,n)} 
&=& {1\over r+1}\left( {n\choose r} - \size{\mathcal{U}(r,n)} - (n-r+1)\size{\mathcal{U}(r-1,n)}\right)
\end{eqnarray*}
with initial values
$\size{\mathcal{U}(3,n)} = {n\choose 3}/(n-2)$ and 
$\size{\mathcal{U}(4,n)} = {n\choose 4}/(n-2)$.
\end{theorem}

\begin{proof}
Let $n=2^m-1$.
Given an integer $i$, let $\binary{i}$ be its binary representation. Given two integers $i,j$, let
$i\oplus j$ be binary addition in $\mathbbm{Z}_2$.
In this proof we assume the set $S_n=\{1,2,\ldots ,n\}$.
Let $\mathcal{U}(r,n) :=\{ \{i_1, \cdots ,i_r\} \subset S_n : \binary{i_1} \oplus \cdots \oplus \binary{i_r} = \binary{0} \}$. 
Since $n=2^m-1$, all possible binary vectors, except the $0$ vector, 
of length $m$ are possible in $\mathcal{U}(r,n)$.
It is clear that if $A,B \in\mathcal{U}(r,n)$ and $A\neq B$, then $|A \cap B| \leq r-2$. 
For otherwise there must exist $i_1, \dots ,i_{r-1},i_r,i_{r+1}$ 
such that $A=\{i_1,\dots ,i_{r-1},i_{r}\}$ and $B=\{i_1,\dots ,i_{r-1},i_{r+1}\}$. 
This gives $\binary{i_1} \oplus \cdots \oplus \binary{i_{r-1}} = \binary{i_r} = \binary{i_{r+1}}$. 
But since $i_r$ and $i_{r+1}$ are distinct, this is a contradiction.

Thus if we take any subset $\mathcal{V} \subseteq \mathcal{U}(r,n)$, 
and then insert into $\mathcal{V}$ all those $r-1$ element subsets of $\sn$ which are 
not contained in any member of $\mathcal{V}$, the resulting collection will be a 
$(r-1)$-partition of $\sn$ and hence represents the hyperplanes of a rank-$r$ paving matroid on $\sn$. 
There are $2^{|\mathcal{U}(r,n)|}$ choices for $\mathcal{V}$, hence we have 
$\size{\Matroids{r-1}{r}{\sn}} \geq 2^{|\mathcal{U}(r,n)|}$. 
It is now required to find $\size{\mathcal{U}(r,n)}$,
\begin{eqnarray*}
\lefteqn{\size{\mathcal{U}(r+1,n)} }\\
&=& \size{\{ \{i_1, \cdots ,i_{r+1}\} \subset S_n : 
	\binary{i_1} \oplus \cdots \oplus \binary{i_{r+1}} = \binary{0} \}} \\
&=& {1\over r+1} \size{\{ \{i_1, \cdots ,i_{r}\} \subset S_n : 
	\binary{i_1} \oplus \cdots \oplus \binary{i_{r}} = \binary{z}, \; z \in \sn-\{i_1,\dots ,i_r\} \}} \\
&=& {1\over r+1}\left( \size{\{ \{i_1, \cdots ,i_{r}\} \subset S_n\}} \right. \\ 
	&& \left. - \size{\{ \{i_1, \cdots ,i_{r}\} \subset S_n : 
	\binary{i_1} \oplus \cdots \oplus \binary{i_{r}} = \binary{z}, \; z \in \{0,i_1,\dots ,i_r\} \}} \right)\\
&=& {1\over r+1}\left( {n\choose r} - \size{\{ \{i_1, \cdots ,i_{r}\} \subset S_n :
	\binary{i_1} \oplus \cdots \oplus \binary{i_{r}} = \binary{0} \}} \right. \\ 
	&& \left. -\size{\{ \{i_1, \cdots ,i_{r}\} \subset S_n : 
	\binary{i_1} \oplus \cdots \oplus \binary{i_{r}} = \binary{z}, \; z \in \{i_1,\dots ,i_r\} \}} \right)\\
&=& {1\over r+1}\left( {n\choose r} - \size{\mathcal{U}(r,n)} \right. \\ 
	&& \left. - (n-(r-1))\size{\{ \{i_1, \cdots ,i_{r-1}\} \subset S_n : 
	\binary{i_1} \oplus \cdots \oplus \binary{i_{r-1}} = \binary{0}\}} \right)\\
&=& {1\over r+1}\left( {n\choose r} - \size{\mathcal{U}(r,n)} - (n-r+1)\size{\mathcal{U}(r-1,n)}\right).
\end{eqnarray*}
The initial conditions for this recursion are
$\size{\mathcal{U}(3,n)} = {n\choose 3}/(n-2)$ and 
$\size{\mathcal{U}(4,n)} = {n\choose 4}/(n-2)$.
$\mathcal{U}(3,n)$ of course corresponds to certain 2-partitions of $\sn$, 
which are the hyperplanes of rank-3 paving matroids.
\end{proof}

\section{Concerning the points-lines-planes conjecture}
The following conjecture is stated in Seymour~\cite{seymour}. Originally Welsh~\cite[p.289]{welshbook} introduced the phrase `Points-lines-planes
conjecture' for the $k=2$ instance of Mason~\cite[Conjecture 1C]{mason}. In the current form it is actually the $k=2$ instance of the slightly more general
Mason~\cite[Conjecture 1D]{mason}.
\begin{conjecture}[Points-Lines-Planes]
For all $M\in \Matroids{2}{4}{\sn}$,
\begin{eqnarray*}
{W_2(M)}^2 &\geq& {3(W_1(M)-1) \over 2(W_1(M)-2)} W_1(M) W_3(M),
\end{eqnarray*}
with equality if and only if $\F_3(M) = {\sn\choose 3}$.
\end{conjecture}
The most general case proven so far is:
\begin{theorem}[Seymour~\cite{seymour}] 
The points-lines-planes conjecture is true for all 
$\{M \in \Matroids{2}{4}{\sn}\;:\; |F|\leq 4 \mbox{ for all } F \in \F_2(M)\}$.
\end{theorem}

Let us also mention that in related work, Kung~\cite{kung} proved 
that if $M$ is a simple matroid of rank at least 5,
in which all lines contain the same number of points, 
then $W_2(M) \leq W_3(M)$.

The purpose of this section is to prove that the points-lines-planes 
conjecture is true for all simple rank-4 matroids if and only if it is true for 
a special subcollection of them. 

Let us briefly outline Knuth's~\cite{knuthrandom} algorithm for constructing all matroids on the finite set $\sn$ 
(i.e. for constructing the class $\cup_{r=1}^{n}\Matroids{1}{r}{\sn}$).
We restate the algorithm in a slightly different way to the original paper.
Let $\T$ be a collection of subsets of $\sn$.
Let $\Expand(\T) := \{A\cup a: A\in\T \mbox{ and }a \in\sn-A\}$.
Let $\Random(\T)$ be a collection of random sets where each random set properly contains some member of $\T$.
Let $\Refine(\T,\T ')$ be the collection of sets $\T$ obtained after applying the following:
\ \\[1em]
\centerline{{\parbox{325pt}{Repeat until $A \cap B \subseteq C$ for some $C \in \T '$ for all $A,B \in \T$:
	if $A,B \in \T$, $A \neq B$ and $A \cap B \not\subseteq C$
        for any $C \in \T '$ then replace $A$ and $B$ in $\T$ by the set $A \cup B$.}}}
\ \\[1em]
\centerline{\fbox{\parbox{325pt}{
{\centerline{\underline{Knuth's random matroids algorithm}}}
\begin{enumerate}
\item Let $r=0$ and $\F_0=\{\emptyset\}$.
\item Let $\A = \Expand(\F_r) \cup \Random(\F_{r})$.
\item Let $\F_{r+1} = \Refine(\A,\F_r)$.
\item If $S_n \in \F_{r+1}$ then let $\F = \F_0 \cup \ldots \cup \F_{r+1}$ and Stop.\\
        Otherwise let $r=r+1$ and go to step (2).
\end{enumerate}
}}}
\ \\
\begin{theorem}[Knuth~\cite{knuthrandom}]\label{changethm} \ \\ \vspace*{-2em}\ \\
\renewcommand{\theenumi}{\roman{enumi}}
\begin{enumerate}
\item The resulting collection $\F$ is the collection of flats of a matroid on $\sn$.
\item The resulting matroid obtained from the algorithm does not depend 
on the order of replacements in $\Refine(\A,\F_r)$.
\item Every matroid in $\cup_{r=1}^n \Matroids{1}{r}{\sn}$ is obtainable using this algorithm.
\end{enumerate}
\renewcommand{\theenumi}{\arabic{enumi}}
\end{theorem}

Let $M$ be a simple rank-$(r+1)$ matroid on $S_n$ with hyperplanes $\F_r$,
and let $\F'_{r+1}= \Refine(\Expand(\F_r) \cup \Random(\F_{r}),\F_r)$.
The matroid $M'$ with flats $\F(M) \cup \F'_{r+1}$ is an {erection} of $M$.
If $M'= M$ (this happens if $\F'_{r+1} = \{S_n\}$) then $M'$ is the {\it{trivial erection}}.
If $\Random(\F_{r}) = \emptyset$, then the resulting matroid $M'$ is 
termed the {\it{free erection}} of $M$, denoted $\Free(M)$.

The free erection was first defined in Crapo~\cite{crapo70i} and an alternative construction of it can be found in Nguyen~\cite{nguyen}.
We prove the following result using the operations defined above. 
An alternative proof may be given using the results and terminology of \cite{crapo70i,nguyen}.

\begin{theorem} \label{myequiv}
The points-lines-planes conjecture is true if and only if it is true for the 
class of rank-4 matroids that are the free erection of some rank-3 matroid.
\end{theorem}

\begin{proof}
If the points-lines-planes conjecture is true, then clearly it is true for 
all those rank-4 matroids that are the free erection of some rank-3 matroid, i.e. $\{ \Free(M): M \in \Matroids{2}{3}{\sn},\,\Free(M)\neq M  \}\subset \Matroids{2}{4}{\sn}$.

Conversely, suppose the points-lines-planes conjecture is true for 
all $\{ \Free(M): M \in \Matroids{2}{3}{\sn},\,\Free(M)\neq M \}$.
Every matroid $N \in \Matroids{2}{4}{\sn}$ is the erection of some $M \in \Matroids{2}{3}{\sn}$, 
that is to say, if $M=U_{3}(N)$ then there exists a collection $\Random(\F_2(M))$ (not necessarily unique) such that 
\begin{eqnarray*}
\F_3(N) &=& \Refine(\Expand(\F_2(M))\cup \Random(\F_2(M)) , \F_2(M)).
\end{eqnarray*}
We now show that for all matroids $M \in \Matroids{2}{3}{\sn}$, 
\begin{eqnarray*}
|\Refine(\Expand(\F_2(M))\cup \Random(\F_2(M)),\F_2(M))| &\leq& |\F_3(\Free(M))|.
\end{eqnarray*}
Let $M \in \Matroids{2}{3}{\sn}$ with $\F_2(M) = \F_2$, so that 
\begin{eqnarray*}
\F_3(N) &:=& \Refine(\Expand(\F_2)\cup \Random(\F_2),\F_2).
\end{eqnarray*}
Consider the collections 
$\Refine(\Expand(\F_2) \cup \Random(\F_2),\F_2)$ and\\
$\Refine(\Expand(\F_2),\F_2)$.
Let $\Expand(\F_2) = \{E_1,\ldots ,E_m\}$ and define $E_i'=E_i$ for all $1\leq i\leq m$.
For each $1\leq i\leq m$, if $X\in\Random(\F_2)$ and $X\supseteq E_i$, then let $E_i' = E_i' \cup X$ and
remove $X$ from $\Random(\F_2)$. 
Once $\Random(\F_2)$ is empty, we have $E_i \subseteq E_i'$ for all $i$. 

Recall from Theorem~\ref{changethm} that if $\A$ is a well defined collection in terms of $\F_2$, then
$\Refine(\A-\{A_1,A_2\}\cup\Refine(\{A_1,A_2\},\F_2),\F_2)=
\Refine(\A-\{A_3,A_4\}\cup\Refine(\{A_3,A_4\},\F_2),\F_2)$ for all $A_1,A_2,A_3,A_4 \in\A$ s.t. $A_1\neq A_2$ and $A_3 \neq A_4$.

For each $E_i,E_j \in \Expand(\F_2)$, if $\Refine(\{E_i,E_j\},\F_2) = \{E_i \cup E_j\}$ then 
$\Refine(\{E_i' \cup E_j'\},\F_2) = \{E_i' \cup E_j'\}$ since $E_i \cap E_j \subseteq E_i' \cap E_j'$.
Hence, $\Refine(\{E_1,\ldots ,E_m\},\F_2)$ has at least as many sets as 
\begin{eqnarray*}
\Refine(\{E_1',\ldots ,E_m'\},\F_2) &=& \Refine(\Expand(\F_2)\cup\Random(\F_2),\F_2)
\end{eqnarray*}
and so 
\begin{eqnarray*}
W_3(N)& = &\size{\Refine(\Expand(\F_2(M))\cup\Random(\F_2(M)),\F_2(M))} \\
&\leq & \size{\Refine(\Expand(\F_2(M)),\F_2(M))}\\
& = & W_3(\Free(M)).
\end{eqnarray*}
Since $W_2(M)=W_2(N)$ and $W_1(M)=W_1(N)$, the result follows.
\end{proof}

Let us also observe that since every plane (rank-3 flat) is determined by two unique lines (rank-2 flats), $W_3(M) \leq W_2(M)^2/2$.
The benefit of Theorem~\ref{myequiv} is apparent by considering the number of simple rank-3 matroids having non-trivial erections relative to
the number of simple rank-4 matroids.
Using Theorem~\ref{slc-bound} we have $\size{\Matroids{2}{4}{\sn}} \geq \size{\Matroids{3}{4}{\sn}} \geq 2^{{n\choose 4}/2n}$
and $\size{\{ \Free(M): M \in \Matroids{2}{3}{\sn},\,\Free(M)\neq M \}}\leq \size{\Matroids{2}{3}{\sn}} \leq 2^{{n\choose 3}}$, 
since all rank-3 matroids are uniquely expressible by their collection of bases, a subset of ${S_n\choose 3}$.
Hence
\begin{eqnarray*}
\dfrac{\size{\{ \Free(M): M \in \Matroids{2}{3}{\sn},\,\Free(M)\neq M \}}}{\size{\Matroids{2}{4}{\sn}}}
& \leq & 2^{-(7n+3)(n-1)(n-2)/48}.
\end{eqnarray*}

\section{Numerical values}
As mentioned in Section 2, the tables for $\size{\NIMatroids{2}{r}{\sn}}$ and 
$\size{\NIMatroids{0}{r}{\sn}}$ were first given in Blackburn, Crapo and 
Higgs~\cite{blackcrapohiggs} and Acketa~\cite{acketa84}, respectively. 
These two tables along with the tables for $\size{\Matroids{0}{r}{\sn}}$, $\size{\Matroids{1}{r}{\sn}}$, 
$\size{\NIMatroids{1}{r}{\sn}}$ and $\size{\Matroids{2}{r}{\sn}}$ were given in Dukes~\cite{dukesphd}. 
A modification to the program used to generate these gives the two new tables for 
the number of paving matroids, $\size{\Matroids{r-1}{r}{\sn}}$, and non-isomorphic 
paving matroids, $\size{\NIMatroids{r-1}{r}{\sn}}$ on at most eight elements.

\begin{tiny}
\begin{center}
\begin{tabular}{r||ccccccccc}
\multicolumn{1}{l}{} & \multicolumn{9}{c}{$\size{\Matroids{0}{r}{\sn}}$} \\ \hline
$r$ $\backslash$ $n$ & 0 & 1 & 2 & 3 & 4 & 5 & 6 & 7 & 8 \\
\hline \hline
0  & 1 & 1 & 1 & 1 & 1  & 1  & 1    & 1     & 1 \\
1  &   & 1 & 3 & 7 & 15 & 31 & 63   & 127   & 255       \\
2  &   &   & 1 & 7 & 36 & 171& 813  & 4012  & 20891     \\
3  &   &   &   & 1 & 15 & 171& 2053 & 33442 & 1022217   \\
4  &   &   &   &   & 1  & 31 & 813  & 33442 & 8520812 \\
5  &   &   &   &   &    & 1  & 63   & 4012  & 1022217   \\
6  &   &   &   &   &    &    & 1    & 127   & 20891     \\
7  &   &   &   &   &    &    &      & 1     & 255       \\
8  &   &   &   &   &    &    &      &       & 1 \\
\hline
$\size{\Matroids{0}{}{\sn}}$& 1 & 2 & 5 & 16 & 68 & 406 & 3807 & 75164 & 10607540  \\
\hline
\end{tabular}
\end{center}

\begin{center}
\begin{tabular}{r||ccccccccc}
\multicolumn{1}{l}{} & \multicolumn{9}{c}{$\size{\NIMatroids{0}{r}{\sn}}$}\\ \hline
$r$ $\backslash$ $n$ & 0 & 1 & 2 & 3 & 4 & 5 & 6 & 7 & 8 \\
\hline \hline
0 &  1& 1 & 1 & 1 & 1 & 1 & 1 & 1 & 1\\
1 &   & 1 & 2 & 3 & 4 & 5 & 6 & 7 & 8\\
2 &   &   & 1 & 3 & 7 & 13& 23& 37& 58\\
3 &   &   &   & 1 & 4 & 13& 38&108&325\\
4 &   &   &   &   & 1 & 5 & 23&108&940\\
5 &   &   &   &   &   & 1 & 6 & 37&325\\
6 &   &   &   &   &   &   & 1 & 7 & 58\\
7 &   &   &   &   &   &   &   & 1 & 8 \\
8 &   &   &   &   &   &   &   &   & 1 \\
\hline
$\size{\NIMatroids{0}{}{\sn}}$ & 1 & 2 & 4 & 8 & 17 & 38 & 98 & 306 & 1724 \\
\hline
\end{tabular}
\end{center}

\begin{center}
\begin{tabular}{r||cccccccc}
\multicolumn{1}{l}{} & \multicolumn{8}{c}{$\size{\Matroids{1}{r}{\sn}}$} \\
\hline
$r$ $\backslash$ $n$ & 1 & 2 & 3 & 4 & 5 & 6 & 7 & 8 \\
\hline \hline
1 &  1 & 1 & 1 &  1 &   1 &    1 &     1 &     1   \\
2 &    & 1 & 4 & 14 &  51 &  202 &   876 &   4139  \\
3 &    &   & 1 & 11 & 106 & 1232 & 22172 & 803583  \\
4 &    &   &   &  1 &  26 &  642 & 28367 & 8274374 \\
5 &    &   &   &    &   1 &   57 &  3592 & 991829  \\
6 &    &   &   &    &     &    1 &   120 &  19903  \\
7 &    &   &   &    &     &      &     1 &    247  \\
8 &    &   &   &    &     &      &       &      1  \\
\hline
$\size{\Matroids{1}{}{\sn}}$&1& 2 & 6 & 27 & 165 & 2135 & 55129 & 10094077 \\
\hline
\end{tabular}
\end{center}

\begin{center}
\begin{tabular}{r||cccccccc}
\multicolumn{1}{l}{} &
\multicolumn{8}{c}{$\size{\NIMatroids{1}{r}{\sn}}$} \\
\hline
$r$ $\backslash$ $n$ & 1 & 2 & 3 & 4 & 5 & 6 & 7 & 8 \\
\hline \hline
1 &  1 & 1 & 1 & 1 & 1 & 1 & 1 & 1  \\
2 &    & 1 & 2 & 4 & 6 &10 &14 & 21 \\
3 &    &   & 1 & 3 & 9 &25 &70 & 217 \\
4 &    &   &   & 1 & 4 &18 &85 & 832  \\
5 &    &   &   &   & 1 &5  &31 & 288  \\
6 &    &   &   &   &   & 1 & 6 & 51 \\
7 &    &   &   &   &   &   & 1 & 7  \\
8 &    &   &   &   &   &   &   & 1  \\
\hline
$\size{\NIMatroids{1}{}{\sn}}$ & 
   1& 2 & 4 & 9 & 21&60 &208&1418\\
\hline
\end{tabular}
\end{center}

\begin{center}
\begin{tabular}{r||ccccccc}
\multicolumn{1}{l}{} & \multicolumn{7}{c}{$\size{\Matroids{2}{r}{\sn}}$} \\ \hline
$r$ $\backslash$ $n$ & 2 & 3 & 4 & 5 & 6 & 7 & 8\\
\hline \hline
2 &    1 & 1 & 1 & 1  & 1   & 1    & 1  \\
3 &      & 1 & 5 & 31 & 352 & 8389 & 433038 \\
4 &      &   & 1 & 16 & 337 & 18700& 7642631 \\
5 &      &   &   & 1  & 42  & 2570 & 907647 \\
6 &      &   &   &    & 1   & 99   & 16865 \\
7 &      &   &   &    &     & 1    & 219 \\
8 &      &   &   &    &     &      & 1 \\
\hline
$\size{\Matroids{2}{}{\sn}}$& 1 & 2& 7 & 49 & 733 & 29760& 9000402 \\
\hline
\end{tabular}
\end{center}

\begin{center}
\begin{tabular}[t]{r||ccccccc}
\multicolumn{1}{l}{} & \multicolumn{7}{c}{$\size{\NIMatroids{2}{r}{\sn}}$} \\ \hline
$r$ $\backslash$ $n$ & 2 & 3 & 4 & 5 & 6 & 7 & 8\\
\hline \hline
2 & 1 & 1 & 1 & 1 & 1 & 1 & 1 \\
3 &   & 1 & 2 & 4 & 9 & 23& 68  \\
4 &   &   & 1 & 3 & 11& 49& 617 \\
5 &   &   &   & 1 & 4 & 22& 217 \\
6 &   &   &   &   & 1 & 5 & 40  \\
7 &   &   &   &   &   & 1 & 6  \\
8 &   &   &   &   &   &   & 1 \\
\hline
$\size{\NIMatroids{2}{}{\sn}}$& 1 & 2 & 4 & 9 & 26 & 101 & 950 \\
\hline
\end{tabular}
\end{center}

\begin{center}
\begin{tabular}{r||ccccccc||ccccccc}
\multicolumn{1}{l}{} & \multicolumn{7}{c}{$\size{\Matroids{r-1}{r}{\sn}}$} & \multicolumn{7}{c}{$\size{\NIMatroids{r-1}{r}{\sn}}$} \\ \hline
$r$ $\backslash$ $n$ & 2 & 3 & 4 & 5 & 6 & 7 & 8 & 2 & 3 & 4 & 5 & 6 & 7 & 8\\
\hline \hline
2 &    1 & 1 & 1 & 1 & 1 & 1 & 1  &
  1 & 1 & 1 & 1 & 1 & 1 & 1 \\
3 &      & 1 & 5 & 31 & 352 & 8389 & 433038 &
    & 1 & 2 & 4 & 9 & 23& 68  \\
4 &      &   & 1 & 6 & 82 & 6149 & 4464328  &
    &   & 1 & 2 & 5 & 18 & 322 \\
5 &      &   &   & 1  & 7  & 239 & 239173 &
    &   &   & 1 & 2 & 5& 39 \\
6 &      &   &   &    & 1   & 8   & 772 &
    &   &   &   & 1 & 2 & 6  \\
7 &      &   &   &    &     & 1    & 9 &
    &   &   &   &   & 1 & 2  \\
8 &      &   &   &    &     &      & 1 &
    &   &   &   &   &   & 1 \\
\hline
Total & 1 & 2& 7 & 39 & 443 & 14787 & 5137322  &
 1 & 2 & 4 & 8 & 18 & 50 & 439 \\
\hline
\end{tabular}
\end{center}
\end{tiny}
\ \\

The program is freely available by e-mailing the author. It is written in the C-programming language and is easily compiled on FreeBSD and Linux platforms. 

\section*{Acknowledgments}
The author would like to thank Francesco Brenti and the anonymous referee for helpful comments, and Christian Krattenthaler for an interesting suggestion.
This work was funded by the Research Training Network `Algebraic Combinatorics in Europe', grant HPRN-CT-2001-00272.

\end{document}